\documentclass{amsart}
\usepackage{ifpdf}
\usepackage{amsmath, amsthm, amsfonts, amssymb, amscd}
\usepackage[active]{srcltx}
\usepackage[pdftex]{graphicx}
    \ifpdf

      \usepackage[pdftex]{graphicx}  

      \usepackage[pdftex]{hyperref}

    \else

      \usepackage[dvips]{graphicx}  

      \newcommand{\href}[2]{#2}

    \fi

\newcommand{\abs}[1]{\left\lvert{#1}\right\rvert}

\DeclareMathOperator{\inter}{\rm{int}}

\DeclareMathOperator{\cl}{cl}

\DeclareMathOperator{\per}{\rm{Per}}

\DeclareMathOperator{\fix}{\rm{Fix}}
\DeclareMathOperator{\ine}{\rm{Ine}}
\DeclareMathOperator{\ess}{\rm{Ess}}

\newcommand{\ol}{\overline}
\renewcommand{\hat}{\widehat}
\newcommand{\til}{\widetilde}

\newcommand{\R}{\mathbb{R}}\newcommand{\N}{\mathbb{N}}
\newcommand{\Z}{\mathbb{Z}}\newcommand{\Q}{\mathbb{Q}}
\newcommand{\T}{\mathbb{T}}
\newcommand{\A}{\mathbb{A}}

\newcommand{\sm}{\setminus}

\newtheorem{theorem}{Theorem}
\newtheorem{corollary}[theorem]{Corollary}
\newtheorem{lemma}[theorem]{Lemma}
\newtheorem{proposition}[theorem]{Proposition}

\newtheorem{question}[theorem]{Question}

\newtheorem*{theorem*}{Theorem}

\theoremstyle{definition}

\theoremstyle{remark}

\begin{document}

\title{ On non-contractible periodic orbits for surface homeomorphisms}
\author{F\'abio Armando Tal}
\address{Instituto de Matem\'atica e Estat\'\i stica, Universidade de S\~ao Paulo,  Rua do Mat\~ao 1010, Cidade Universit\'aria, 05508-090 S\~ao Paulo, SP, Brazil}
\email{fabiotal@ime.usp.br}
\thanks{Supported by CNPq grant 304474/2011-8 and Fapesp 2011/16265-8}

\begin{abstract}
In this work we study homeomorphisms of closed orientable surfaces homotopic to the identity, focusing on the existence of non-contractible periodic orbits.  We show that, if $g$ is such a homeomorphism, and if $\hat g$ is its lift to the universal covering of $S$ that commutes with the deck transformations, then one of the following three conditions must be satisfied: (1) The set of fixed points for $\hat g$ projects to a closed subset $F$ which contains an essential continuum, (2) $g$ has non-contratible periodic points of every sufficiently large period, or (3) there exists an uniform bound $M>0$ such that, if $\hat x$ projects to a contractible periodic point, then the $\hat g$ orbit of $\hat x$ has diameter less or equal to $M$. Some consequences for homeomorphisms of surfaces whose rotation set is a singleton are derived.
\end{abstract}

\maketitle

\section{Introduction}

Let $M$ be a compact manifold and $I_h(x,t):M\times [0,1]\to M$ be an isotopy between the identity $Id(x)=I_h(x,0)$ and a map $h(x)=I_h(x,1)$. A fixed point $p$ for $h$ is called contractible if the loop $\gamma_p^h(t)=I_h(p, t),\, t\in [0,1]$ is homotopic, in $M$, to a point.

The study of the existence of contractible fixed points and contractible periodic orbits for homeomorphisms of manifolds has a long history in mathematics. In the particular case that $M$ is a sympletic manifold, $h$ is the time one map of a $\mathcal{C}^2$ periodic Hamiltonian flow and $h$ has finitely many fixed points, the Conley conjecture on the existence of contractible periodic orbits of arbitrary large period was a great motivation for the fields of topological dynamics on surfaces (\cite{Fr92, FH, patrice, patriceduke}) and sympletic dynamics (see \cite{Hi, Gi} for some recent results).

Lately, the study of the existence of \emph{non-contractible} periodic points has been the focus of significant
work in sympletic geometry (\cite{Gal, BPS, Joa1, Gurel, Pedrojoa}) with great developments, usually employing Floer homology techniques. This note is concerned  with the existence of such points for general homeomorphisms of surfaces. Our intention is to understand the dynamical restrictions that are implied whenever all $h$--periodic points are contractible. Our main result  is:
\begin{theorem}\label{th:pontos periodicos}
Let $S$ be a compact orientable surface, $g:S\to S$ a homeomorphism isotopic to the identity and $\hat g$ a lift of $g$ to the universal covering $\hat S$ that commutes with the deck transformations. Then either:
\begin{enumerate}
\item{ $\pi(\fix(\hat g))$ is essential,}
\item{there exists $T\not= Id$ a deck transformation such that, for every $\frac{p}{q}\in\Q$ with $\abs{\frac{p}{q}}$ sufficiently small, there exists $\hat z=\hat z(\frac{p}{q})$ such that $\hat g^q(\hat z)=T^p(\hat z)$, or}
\item{there exists $M>0$ such that, if $ \hat x$ is $\hat g$ periodic, then the orbit of $\hat x$  has diameter less then $M$.} 
\end{enumerate}
\end{theorem}
and a direct consequence is:

\begin{corollary}
Let $S$ be a compact orientable surface, $g:S\to S$ a homeomorphism isotopic to the identity with finitely many fixed points, and $\hat g$ a lift of $g$ to the universal covering $\hat S$ that commutes with the deck transformations. Then either:
\begin{enumerate}
\item {there exists $M>0$ such that, if $ \hat x$ is $\hat g$ periodic, then the orbit of $\hat x$  has diameter less then $M$ or,}
\item{ $g$ has non--contractible periodic orbits of arbitrary large period.}
\end{enumerate}
\end{corollary}

The previous theorem has some other interesting consequences. Whenever a homeomorphism $f$ of a compact orientable surface $S$ leaves invariant a measure of full support $\mu$ (in which case the set of nonwandering points of $f,\, \Omega(f)$, is the whole surface), the prevailing picture of the dynamics is that of a phase space partitioned in two different regions: one of points contained in periodic topological disks, often denoted as ``elliptic islands'', and another connected compact set where the dynamics is transitive. Recent works (\cite{Tobias, Inventiones, Salvador}) have shown that for torus homeomorphisms this description is often valid. A topological description of the ``islands'' involves the notion of \emph{inessential} points for the dynamics, those which are either wandering or contained in periodic topological disks. Denote by $\ine(f)$ the set of inessential points, and by $\Omega(f)$ the set of nonwandering points. Inessential nonwandering points in many ways behave as periodic points, and the conclusions of Theorem \ref{th:pontos periodicos} can be adapted to them:

\begin{theorem}\label{th:pontos inessenciais}
Let $S$ be a compact orientable surface, $g:S\to S$ a homeomorphism isotopic to the identity and $\hat g$ a lift of $g$ to the universal covering $\hat S$ that commutes with the deck transformations. Then either:
\begin{enumerate}
\item{ $\pi(\fix(\hat g))$ is essential,}
\item{$g$ has non--contractible periodic orbits,}
\item{there exists $M>0$ such that, if $x\in \ine(g)\cap \Omega(g)$ and $\pi(\hat x)=x$, then the orbit of $\hat x$  has diameter less than $M$.} 
\end{enumerate}
\end{theorem}

Another interesting consequence is obtained for \emph{irrotational} homeomorphisms of compact surfaces, i.e., those  homeomorphisms homotopic to the identity for which every point have null rotation vector (a precise definition of irrotacional homeomorphism is given in the next section). These homeomorphisms are in some sense the simplest possible for the rotational viewpoint, and for them all periodic orbits must be contractible. Irrotational homeomorphisms that are also area preserving are necessarily Hamiltonian in the sense of \cite{patriceduke}, and  the dynamics may be very complicated when the set of fixed points is large (see \cite{exemplo}). Apart from this case, a reasonable question is:

\begin{question}
If $S$ is a compact surface, $g$ is irrotational and preserves a measure  of full support, and $\fix(g)$ is contained in a topological disk, is there $M>0$ such that all $\hat g$ orbits have diameter less then $M$?
\end{question}
Theorem \ref{th:pontos periodicos} shows this to be true if periodic points are dense in $S$. When $S$ is the torus, the main result of  \cite{pseudorotacoes}, Theorem \ref{th:pseudorotacoes}, indicates this is probably true, and Theorem \ref{th:pontos inessenciais} allows us to improve it:

\begin{theorem}\label{th:irrotacional}
Let $g:\T^2\to\T^2$ be an irrotational homeomorphisms preserving a measure of full support, and $\hat g$ its lift to $\R^2$. If $\fix(g)$ is contained in a topological disk, then there exists $v\in\Z^2$ such that the diameter of every orbit by $\hat g$ is uniformly bounded in the $v$ direction.
\end{theorem}

The next section introduces notation and the main background needed for the proof of Theorem \ref{th:pontos periodicos}, Nielsen-Thurston theory and Le Calvez work on equivariant Brouwer lines, and the  proofs are done in the last section.

\section{Preliminar results}
In what follows $S$ will be a compact orientable surface, possibly with boundary, $\hat S$ its universal covering and $\pi:\hat S\to S$ the covering projection.  We denote the set of deck transformations of $\hat S$ by $Deck(S)$.  We endow $S$ with a metric $d_S$ and denote the induced metric in $\hat S$ by $d_{\hat S}$.

Following \cite{Inventiones}, an open set $U\subset S$ is called \emph{inessential} if every closed loop in $U$ is homotopically trivial in $S$. A set $F\subset S$ is called inessential if it has an inessential neighborhood, otherwise $F$ is called essential.  A connected set $F\subset S$ is called bounded if whenever $\hat F$ is a connected components of $\pi^{-1}(F)$, then $\hat F$ is bounded in $\hat S$, in which case the diameter of $F$ is defined as the diameter of $\hat F$.  While open inessential sets need not be bounded, it is immediate that any \emph{closed} inessential set is contained in a closed bounded topological disk $\overline D$, and as such the diameter of all its connected components is uniformly bounded by the diameter of $\overline D$.

Given a homeomorphism $h: S\to S$, we say that $x\in S$ is \emph{inessential} for $h$ if there exists an $h$-invariant neighborhood of $x$ which is inessential, otherwise $x$ is called essential. The set of essential and inessential points of $h$ are denoted by $\ess(h)$ and $\ine(h)$, respectively. If $h$ is homotopic to the identity, and $\hat h$ is a fixed lift of $h$ that commutes with deck transformations, we say that the orbit of $x$ is bounded if, for any $\hat x\in\pi^{-1}(x)$, the set $\bigcup_{i\in\Z}\{\hat h^{i}(\hat x)\}$ is bounded. We denote the set of nonwandering points of $h$ by $\Omega(h)$.

Whenever $h$ is isotopic to the identity, we chose an isotopy $I_h:S\times t\to S$ satisfying $I_h(x,0)=x$ and $I_h(x,1)=h(x)$. We denote by $I_{\hat S}$ the induced isotopy on the covering space, and by $\hat h$ the corresponding lift of $h$.  We also denote by $\gamma_x^h:[0,1]\to S, \gamma_x^h(t)=I_h(x,t)$ and we write just $\gamma_x$ whenever there is no possibility for confusion.

Given curves $\alpha,\beta:[0,1]\to S$ we denote their respective images by $[\alpha],[\beta]$, and their concatenation by $\alpha*\beta$. If $\hat S$ is homeomorphic to $\R^2$, $\alpha:[0,1]\to \hat S$ is a closed curve and $x$ does not belong to $[\alpha]$, we denote by $W(\alpha, x)$ the winding number of $\alpha$ around $x$. Whenever $\hat h$ is a homeomorphism of $\hat S$, $\hat x$ is a periodic point for $\hat h$ with minimal period $n$, $\alpha=\gamma_{\hat x}*\gamma_{\hat h(\hat x)}*...*\gamma_{\hat h^{n-1}(\hat x)}$ and  $\hat p\notin[\alpha]$, we define the index of  the orbit of $\hat x$ with respect to $\hat p$ as $W(\alpha,\hat p)$.

Given $v\in \R^2_*$, we denote by $P_v:\R^2\to \R, P_v(\hat x)= \langle x; v\rangle$, the projection in the $v$ direction. If $S=\T^2$ and $\hat S=\R^2$, we say that the orbit of $x\in\T^2$ is bounded in the $v$ direction if for any $\hat x\in\pi^{-1}(x)$, the set $P_v(\bigcup_{i\in\Z}\{\hat h^{i}(\hat x)\})$ is bounded.

\subsection{Rotation set  for homeomorphisms of the annulus}

Let $\ol \A= S^1\times [0,1]$ be the closed annulus, $f: \ol \A\to\ol \A$ be a homeomorphism preserving orientation and boundary components, and $\widetilde f$ a lift of $f$ to the universal covering $\hat{\ol{\A}}= \R\times[0,1]$, and let $\pi:\hat{\ol{\A}}\to \ol \A$ be the covering map . Given a point $x\in\ol{\A}$, we define, whenever the limit exists, the rotation number of $x$ as $\rho(x,\hat f)=\lim_{n\to\infty}\frac{(\hat f^{n}(\hat x)-\hat x)_1}{n}$, where $\hat x\in\pi^{-1}(x)$ and $(\hat a)_1$ denotes the first coordinate of $\hat a$. Unlike the rotation number of an orientation-preserving circle homeomorphism, the rotation number of a point $x$ does not need to exist. Furthermore, although all points have the same rotation number for a circle homeomorphism, for annulus homeomorphisms different points may have different rotation numbers. Therefore it is useful to define the pointwise rotation set of $\hat f$ as $\rho(\hat f)=\bigcup_{x\in\ol{A}}\{\rho(x,\hat f)\}$, the set of all different rotation numbers. 

The main result in \cite{handelanel} shows the strong relation between rotation numbers and periodic orbits:
\begin{theorem}\label{th:handelanel}
For every rational $\frac{p}{q}\in \rho(\hat f)$ there exists $\hat x$ in $\hat{\ol{A}}$ such that $\hat f^q(\hat x)=\hat x+(p,0)$
\end{theorem}

\subsection{Brouwer--Le Calvez foliation}

The equivariant Brouwer theory developed in \cite{patrice} is a very useful tool in dealing with homeomorphisms of surfaces homotopic to the identity, and we briefly describe it in the context we require.

For our purposes, let  $h$ be a homeomorphism of $S$ isotopic to the identity, $I_h$ an isotopy between $Id_S$ and $h$,  and $I_{\hat h}$ the lifted isotopy between $Id_{\hat S}$ and $\hat h$, such that $\hat h$ commutes with any $T\in Deck(S)$. We assume that $\fix(\hat h)$ is totally disconnected, and that there exists a $X\subset \pi(\fix(\hat h))$ such that for every $p$ in $X,\, t\in [0,1],\, I_h(p,t)=p$. Let $\mathcal{W}$ be an oriented topological foliation of $S\setminus X$, which we regard as a flow with singularities. We say that $\mathcal{W}$ is dynamically transverse to $I_h$ if, for every $x\in S\setminus X$, the curve $\gamma_x^h$ is homotopic with fixed end points, in $S\setminus X$, to an arc $\beta$ positively transversal to $\mathcal{W}$, that is, that locally crosses every $\mathcal{W}$ leaf from left to right. In this case, if $\hat X=\pi^{-1}(X)$ and $\hat{\mathcal{W}}$ is the lifted foliation $\hat S\setminus \hat X$, then $\hat{ \mathcal{W}}$ is dynamically transverse to $I_{\hat h}$. We call $\mathcal{W}$ a Brouwer--Le Calvez foliation for $(I_h,S\setminus X)$.  The result we use, a consequence of  \cite{patrice} and \cite{jaulent}, is:

\begin{theorem}\label{th:folheacao}
If $h:S\to S$ is homotopic to the identity and $\pi(\fix(\hat h))$ is totally disconnected, then there exists a compact set $X\subset S$, an isotopy $I_h$ between $Id_S$ and $h$ which pointwise fixes $X$ and a Brouwer--Le Calvez foliation $\mathcal{W}$ for $(I_h,S\setminus X)$.
\end{theorem}

The following results follow from the ideas in \cite{abertos}, and are included for completeness:

\begin{proposition}\label{pr:bolaefolhas}
If $\mathcal{W}$ is a Brouwer--Le Calvez foliation for $(I_h, S\setminus X)$, then for any $x\in S\setminus X$ there exists $\varepsilon>0$ such that, if $\mathcal{W}_x$ is the leaf of $\mathcal{W}$ passing through $x$, then:
\begin{itemize}
\item{For any $y\in B_{\varepsilon}(x)$, the curve  $\gamma_{h^{-1}(y)}^h*\gamma_y^h$ is homotopic, with fixed endpoints, to an arc positively transversal to the foliation that crosses $\mathcal{W}_x$ at least once,}
\item{If the restriction of $g:S\to S$ to $S\setminus B_{\varepsilon}(x)$ is equal to $h$, then there exist $I_g$ such that $\mathcal{W}$ is a Brouwer--Le Calvez foliation for $(I_g,S\setminus X)$.}
\end{itemize}
\end{proposition}
\proof
If $M=S\setminus X$ and $\hat M$ is its universal covering space and $\hat h$ the lift of $h$ to $\hat M$, then $\hat h$ has no fixed points. Let $\widehat{\mathcal{W}}$ be the lift of $\mathcal{W}$ to $\hat M$. Since $\widehat{\mathcal{W}}$ is a oriented foliation without singularities of $\hat M$ and $\hat M$ is homeomorphic to the plane, every  leaf  $\widehat\phi$ of $\widehat{ \mathcal{W}}$ is a line, that is, the image of a proper and injective map from $\R$ into $\hat M$, and naturally divides $\hat M$ into two regions, the left of $\widehat \phi$, denoted by $L(\widehat \phi)$ and the right of $\widehat \phi$, denoted by $R(\widehat \phi)$. Since $\widehat{\mathcal{W}}$ is dynamically transverse to $I_{\hat h}$, the pre-image of every leaf is contained in its right, and the image of any leaf must be contained in its left. Therefore every leaf of $\mathcal{W}$ lifts to a Brouwer line in $\hat M$. In particular, if $\hat x$ is a pre-image of $x$ in $\hat M$ and $\hat{W}_{\hat x}$ is leaf of $\hat{\mathcal{W}}$ that passes through $\hat x$, there exists a neighborhood $\hat V$ of $\hat x$ such that $\hat h^{-1}(\hat V)$ and $\hat h(\hat V)$ are on opposite sides of $\hat{W}_{\hat x}$. Both claims follow trivially by choosing $\varepsilon>0$ such that $B_{\varepsilon}(\hat x)\subset \hat V$
\endproof

\begin{proposition}
Let $\Phi:\R\times\R^2$ be a topological flow, $\beta:[0,1]\to \R^2$ a closed curve positively transversal to the flow lines of $\Phi$ with finitely many self-intersections, and $\phi_x$ a flow line intersecting  the image of $\beta$. If $\alpha(\phi_x)$ and $\omega(\phi_x)$ are the $\alpha$ and $\omega$ limits of  the flow line, then at least one of them is nonempty, compact, and contained in a single bounded connected component of $\R^2\setminus  [\beta]$. Furthermore, there exists an integer $l\not=0$ such that either $W(\beta, y)=l$ for every $y$ in $\alpha(\phi_x)$,   or $W(\beta, y)=l$ for every $y$ in $\omega(\phi_x)$. 
\end{proposition}

\proof

Let $\beta_0=\beta$ and define inductively, for $i\ge 1$,
$$ \,t_i=\min\{t\le 1\mid {\hbox{there exists }} 0\le s<t,\, \beta_{i-1}(s)=\beta_{i-1}(t)\},$$ 
and $s_i$ as the point in $[0,t_i)$ such that $\beta_{i-1}(s_i)=\beta_{i-1}(t_i)$. Note that, by the definition of $t_i$ as minimal, $s_i$ is uniquely defined. Let $\widetilde \beta_{i}:[0,1]\to\R^2$ be a positive reparametrization of the $\beta_{i-1}|_{[s_i, t_i]}$, and, provided $t_i\not=1,\,\beta_i:[0,1]\to\R^2$ a positive reparametrization of the concatenation $\beta_{i-1}|_{[0,s_i]}*\beta_{i-1}|_{[t_1,1]}$.  Note also that, since $\beta_0$ had finitely many self--intersections, then there exists $n$ such that $t_{n}=1$ and therefore $\beta_{n-1}=\til\beta_n$. We remark that each $\tilde \beta_i$ is a simple closed curve.
By construction $[\beta_{i-1}]=[\beta_i]\cup[\til\beta_i]$ and, if $z\notin [\beta_{i-1}]$, then 
$W(\beta_{i-1},z)=W(\beta_i,z)+W(\til\beta_i,z)$. In particular, if $z\notin[\beta]$, then $W(\beta,z)=\sum_{i=1}^n W(\til\beta_i,z)$.

But if $\gamma$ is a simple closed curve positively transversal to a flow, then $[\gamma]$ intersects a flow line $\phi_x$ at most once, in which case either $\alpha(\phi_x)$ is a nonempty compact set contained in the interior of $\gamma$, $\omega(\phi_x)$ is contained in the exterior of $\gamma$, and for every $z_1\in\alpha(\phi_x), z_2\in\omega(\phi_x),\, W(\gamma,z_1)=-1$ and $ W(\gamma,z_2)=0$, or $\omega(\phi_x)$ is a nonempty compact and contained in the interior of $\gamma$, $\alpha(\phi_x)$ is contained in the exterior of $\gamma$, and for every $z_1\in\alpha(\phi_x), z_2\in\omega(\phi_x),\, W(\gamma,z_1)=0$ and $ W(\gamma,z_2)=1$. Furthermore, if $[\gamma]$ is disjoint from $\phi_x$, and if $z_1\in\alpha(\phi_x)$ and $z_2\in\omega(\phi_x)$, then $W(\gamma,z_1)=W(\gamma,z_2)$. This implies that, for $1\le i\le n$, $W(\til\beta_i,z_1)\le W(\til\beta_i,z_2)$, with the equality holding only if $[\til \beta_i]$ is disjoint from $\phi_x$. As $[\beta]$ is not disjoint from $\phi_x$, there exists $i_0$ such that $[\til\beta_{i_0}]$ is not disjoint from $\phi_x$ and  $W(\til\beta_{i_0},z_1)< W(\til\beta_{i_0},z_2)$, and so $W(\beta,z_1)=\sum_{i=1}^n W(\til\beta_i,z_1)<\sum_{i=1}^n W(\til\beta_i,z_2)=W(\beta,z_2)$
\endproof

\subsection{Nielsen Thurston Theory}

We need to make use of the Thurston's classification theorem for homeomorphisms of surfaces. We describe only the elements needed for our result, the reader is referred to \cite{flp} for a more comprehensive description. Let $S$ be a compact surface and $X=\{x_1,...,x_n\}$ a collection of a finite number of points $S$ such that the Euler characteristic of $S\setminus X$ is negative. Let $f:S\to S$ be a homeomorphism such that $f(X)=X$ and let $\til S= S\setminus X$. Then $f$ is said to be:
\begin{itemize}
\item{\emph{reducible relative to $X$} if there exists a finite set of disjoint simple closed curves $\{\gamma_1, ...,\gamma_n\}$ such that each $\gamma_i\subset \til S$ is neither null homotopic in $\til S$ nor homotopic to a boundary component or to a topological end of $\til S$, and such that $f$ leaves invariant the set $\bigcup_{i=1}^{n}\gamma_i$.}
\item{\emph{pseudo--Anosov relative to $X$} if there exists a pair of transverse  measured foliations $\mathcal{F}^u,\, \mathcal{F}^s$ of $S$ with finitely many singularities, and a real $\lambda>1$ such that $f$ preserves the foliations and such that if $\mu_u,\, \mu_s$ are the respective transverse measures, then $f_*(\mu_u)=\frac{1}{\lambda}\mu_u$ and $f_*(\mu_s)=\lambda\mu_s$.  Furthermore, every singularity $p$ of this pair of foliations is a $q$-pronged saddle, and if $p\notin X$, then $q\ge 3$.}
\end{itemize}

Thurston's classification theorem  states that, in this setting, for any given homeomorphism $h: S \to S$ such that $h(X)=X$, there exists a homeomorphism $f$ and an  isotopy $H$ from $h$ to $f$ leaving $X$ invariant and such that $f$ is either periodic, reducible relative to $X$ or pseudo--Anosov relative to $X$. 

Whenever $h$ is homotopic to a pseudo--Anosov homeomorphism $f$, then each $f$ orbit must be globally shadowed by an $h$ orbit, and periodic orbits of $f$ are globally shadowed by periodic orbits of $h$. In particular, and most relevant for us, we have the following result, a direct consequence of the main theorem of \cite{handelshadowing}:
\begin{theorem}\label{th:handelshadowing}
Let $f:S\to S$ be pseudo--Anosov relative to $X$,  and let $h:S\to S$ be homotopic to $f$. Let $\widetilde S$ the universal covering of $S\setminus X$ and let $\tilde f$ and $\tilde h$ be homotopic lifts of $f$ and $g$ to $\widetilde S$, respectively. Then, if $x\in \widetilde S,\,T:\widetilde S\to\widetilde S$ is a deck transformation and $n$ are such that $\tilde f^n(x)=T(x)$, then there exists $y\in \widetilde S$ such that $\tilde h^{n}(y)=T(y)$.
\end{theorem}

For homeomorphisms of the closed annulus $\ol \A=S^1\times [0,1]$, we use the following theorem  from \cite{toby} relating Nielsen-Thurston theory and rotation sets :
\begin{theorem}\label{th:toby}
Let $f:\ol \A\to \ol \A$ be a homeomorphism of the annulus pseudo--Anosov relative to an invariant finite set $X\subset \ol \A$, and let $\hat f$ be its lift to the universal covering of $\ol \A$. If, for every $x$ in $X$, we have $\rho(x, \hat f)= 0$, then $0\in \inter(\rho(\hat f))$.
\end{theorem}

Finally, 

\begin{lemma}\label{lm:toby2}
Let $\A=S^1\times \R$ be the \emph{open} annulus, $X\subset \A$ be a finite set, $h:\A\to \A$ be a homeomorphisms such that $h(X)=(X)$ and $h$ is homotopic in $\A\setminus X$ to $f$, a homeomorphisms pseudo--Anosov relative to $X$, and  let $\hat h$ be the lift of $h$ to the universal covering $\hat \A$ of $\A$. If for every $x$ in $X,\,\rho(\hat h, x)= 0$, then there exists $c>0$ such that, for any $p, q\in \Z$ with $q\not=0$ and $\abs{\frac{p}{q}}<c$, there exists $\hat y$ such that $\hat h^{q}(\hat y)=(\hat y)+ (p, 0)$.
\end{lemma}
\proof
First note that, from Brouwer theory, since $\hat h$ has periodic points it must also have fixed points, and so the result follows when $p=0$. 

Let $\hat f$ be the lift of $f$ to $\hat \A$ such that, if $x\in X$, then $\rho(\hat f, x)=0$. If $\overline \A$ is the compactification of $A$ by adding two boundary circles at each end of the annulus, then there exists an extension $\ol f$  of $f$ to $\ol \A$ such that $\ol f$ is pseudo--Anosov relative to $X$ in $\ol \A$. Furthermore, there exists a lift $\hat{\ol f}$ of $\ol f$ to $\hat{\ol \A}$ such that for all $x\in \A$, if $\rho(x, \hat f)$ exists, then $\rho(x, \hat{\ol f})$ also exists and is equal to $\rho(x, \hat f)$.

By Theorem \ref{th:toby}, the origin lies in the interior of $\rho(\hat{\ol f})$. Let $c>0$ be such that $[-c,c]\subset \rho(\hat{\ol f})$. Since the rotation number of all points in a boundary component of $\ol{\A}$ is the same, by eventually decreasing $c$, we may assume that the rotation number of the boundary points do not belong to $[-c, 0) \cup (0, c]$. By theorem \ref{th:handelanel}, for every $p, q\in \Z$ with $p\not=0$ and $q\not=0$ and $\abs{\frac{p}{q}}<c$, there exists $\hat z=\hat z(\frac{p}{q})$ such that $\hat{\ol f}^{q}(\hat z)=(\hat z)+ (p, 0)$. By the choice of $c,\, \pi(\hat z)$ must belong to the interior of $\ol{\A}$ and so, by the construction of $\ol{f}$ it follows that $\hat{f}^{q}(\hat z)=(\hat z)+ (p, 0)$. Finally, by the Handel's shadowing theorem \ref{th:handelshadowing}, for every $\hat z (\frac{p}{q})$ there exists $\hat y= \hat y(\frac{p}{q})$ such that $\hat h^{q}(\hat y)=(\hat y)+ (p, 0)$. 
\endproof

\subsection{Strictly toral homeomorphisms and irrotational homeomorphisms of $\T^2$}

Following \cite{Inventiones} a homeomorphisms $h$ of $\T^2$ homotopic to the identity is called annular if there exists $M>0, \,v\in \Z^2,\, k\in \N $ and a lift $\hat h$ of $h$ such that,  for any $\hat x \in \R^2, n\in \Z,\, \abs{P_v(\hat h^n(\hat x)-\hat x)}<M$.   A homeomorphisms $h$ of $\T^2$ is said to be strictly toral if, for every $k\in\N$, $h^{k}$ is not annular and $\fix(h^k)$ is inessential.  The following result follows from Theorems C and D of \cite{Inventiones}:
\begin{theorem}\label{th:strictlytoral}
If $f:\T^2\to\T^2$ is strictly toral and $\Omega(f)=\T^2$, then $\ine(f)$ is an inessential set, and for any open sets $U, V$ intersecting $\ess(f)$ there exists $n>0$ such that $f^{n}(U)\cap V\not=\emptyset$.  
\end{theorem}

A homeomorphisms $h$ isotopic to the identity of a surface $S$ is called \emph{sympletic} if it preserves a Borel probability measure $\mu$ of full support and it is called \emph{hamiltonian} if the $h$-rotation vector of $\mu$ is null (see \cite{patriceduke} for a precise definition). It is called \emph{irrotational} if, for every $\hat x$ in $\hat S,\, \lim_{n\to\infty}\frac{d_{\hat S}(\hat h^n(\hat x),\hat x)}{n}=0$. The main theorem of \cite{pseudorotacoes} gives a description of homeomorphisms of $\T^2$ which are area preserving and irrotational (and, in particular, hamiltonian):

\begin{theorem}\label{th:pseudorotacoes}
If $h:\T^2\to \T^2$ is irrotational and sympletic, then either $h$ is \emph{not} strictly toral, or the $\hat h$ orbit of any point $\hat x$ is bounded. 
\end{theorem}

\section{Proof of the main results}

Assume that neither (1) nor (3) holds, so that $F=\pi(\fix(\hat g))$ is inessential and that there exists a sequence of $\hat g$ periodic points such that the diameter of the orbits is not uniformly bounded. We will show that condition (2) holds.

The following two propositions are direct consequences of  $F$ being closed and  inessential.

\begin{proposition}\label{pr:praux1}
There exists a topological disk $D$ containing $F$ such that every connected component of $\pi^{-1}(D)$ is bounded.
\end{proposition}

\begin{proposition}\label{cai no disco}
There exists open neighborhood $V\subset D$ of  $F$ such that, if $\hat D$ is a connected component of $\pi^{-1}(D)$, $\hat V$ is the connected component of $\pi^{-1}(V)\subset \hat D$ and $\hat x \in \hat V$, then $\hat g (\hat x) \in \hat D$.
\end{proposition}

\begin{proposition}
We may assume $F$ is totally disconnected.
\end{proposition}
\proof
 Since $F$ is contained in $D$ and $D$ is a bounded topological disk, $S\setminus D$ is connected. Let $O$ be the connected component of $S\setminus F$ that contains $S\setminus D$, and let $K$ be $S\setminus O$ ($K$ is the filling of the set $F$). Note that, by a theorem of Brown--Kister (\cite{brown}), since $\pi^{-1}(F)=\fix(\hat g)$, every connected component of $\hat S\setminus \fix(\hat g)$ is invariant. Since the connected components of $\pi^{-1}(K)$ are invariant and uniformly bounded, there exists a $M_1>0$ such that, if the orbit of $\hat x$ has diameter larger than $M_1$, then $\hat x\in \pi^{-1}(O)$.

If $\mathcal{P}$ is the partition of $S$ into sets of the form $\{x\}$, if $x\in O$, and the connected components of $K$, then $\mathcal{P}$ is a upper semicontinuous decomposition of $S$ (see, for instance,  proposition 1.6 of \cite{Inventiones}), and there exists a continuous surjection $\psi:S\to S$, and a homeomorphisms $g':S\to S$ homotopic to the id such that:
\begin{itemize}
\item{$\psi$ is homotopic to the identity,}
\item{$\psi(K)=K'$  is a totally disconnected set,}
\item{$\psi g= g' \psi$ and $\psi\mid_{S\setminus K}$ is a homeomorphisms between $S\setminus K$ and $S\setminus K'$.}
\end{itemize}

Let $\hat{g'}$ be the lift of $g'$ such that $\fix(\hat{g'})=K'$. Then there exists $\hat y\in \hat S$ is such that $(\hat g')^k(\hat y)=T(\hat y)$ if, and only if, there exists $\hat x\in \hat S$ such that $\hat g^k(\hat x)=T(\hat x)$. Furthermore, if $(\hat x_k)_{k\in\N}$ is a sequence of $\hat g$ periodic points  such that the diameter of their orbits is not uniformly bounded, then the same is true for the orbits of the corresponding sequence $\psi(x_k)$ of $\hat {g'}$--periodic points. Therefore conditions (1) and (3) also do not hold for $g'$, and if $g'$ satisfies (2), then so does $g$\endproof

The previous proposition, together with Theorem \ref{th:folheacao}, imply that there exists an isotopy $I_g$, a closed subset $X\subset \pi(\fix(\hat f))$ and a Brouwer--Le Calvez foliation $\mathcal{W}$ for $(I_g, S\setminus X)$.

Let now $D$ and $V$ be sets given by Propositions \ref{pr:praux1} and \ref{cai no disco} and consider the compact set $S\setminus V$. For each $x \in S\sm V$, there exists $\varepsilon(x)>0$, given by proposition \ref{pr:bolaefolhas}, such that, if $y\in B_{\varepsilon(x)}(x)$, then $\gamma_{g^{-1}(y)}*\gamma_y$ is homotopic with fixed points in $S\setminus X$ to an arc transversal to the foliation crossing the leaf $\mathcal{W}_x$. Let $x_1, ..., x_{\ol n}$ and $\varepsilon_i=\varepsilon(x_i), 1\le i\le \ol n$, be such that  $B_{\varepsilon_1}(x_1),...,B_{\varepsilon_{\ol n}}(x_{\ol n})$ is a covering of $S\sm V$ and choose, for each $1\le i\le \ol{n},\,\hat x_i\in \pi^{-1}(x_i)$.

\begin{proposition}\label{pr:visita3}
There exists $M>0$ such that, if the orbit of $\hat x$  has diameter greater then $M$, then there exists $n_1, n_2, n_3\in\Z,\, i\in\{1,...,n\}$ and $T_1, T_2, T_3$ three distinct transformations in $Deck(S)$ such that $\hat g^{n_j}(\hat x)\in T_j(B_{\varepsilon_i}(\hat x_i)), \,j\in\{1,2,3\}$.
\end{proposition}
\proof
 Since $\hat g$ is isotopic to the identity, there exists a constant $C>0$ such that, if the orbit of a point $\hat x$ has diameter greater than $Ck$, then there exists at least $k$ distinct integers $n_1, n_2,..., n_k$ such that $d_{\hat S}(\hat g^{n_s}(\hat x),\hat g^{n_l}(\hat x))>1$ whenever $s\not=l$. Furthermore, by proposition  \ref{cai no disco}, we can assume all of these points lie in $\pi^{-1}(S\setminus V)$. The result follows from taking $M> C(2 \ol{n}+1)$
\endproof

\begin{lemma}\label{lm:enlaca_p_e_tp}
There exists $M>0$ such that, if the orbit of $\hat x\in \per(\hat g)$  has diameter greater then $M$, then there exists $T \in Deck(S),\, T\not= Id,$ and $\hat p$ a singularity of the foliation $\hat{\mathcal{W}}$  such that  the orbit of $\hat x$ has nonzero index with both  $\hat p$ and $T(\hat p)$
\end{lemma}

\proof
By the previous proposition, if the orbit of $\hat x$ has a sufficiently large diameter, then there exists $\hat x_i,\, T_1,T_2$ and $T_3$ distinct deck transformations and $n_1, n_2, n_3$ such that $\hat g^{n_l}(\hat x)\in B_{\varepsilon_i}(T_l(\hat x_i)), \, n_l\in\{n_1, n_2, n_3\}$, and we assume that $0\le n_1<n_2<n_3< k$ where $\hat g^{k}(\hat x)=\hat x$. If $\beta$ is the concatenation 
$\gamma_{\hat x}^{\hat g}*\gamma_{\hat g(\hat x)}^{\hat g} *...*\gamma_{\hat g^{k-1}(\hat x)}^{\hat g}$, then $\beta$ is a closed loop passing through the whole orbit of $\hat x$, and $\beta$ is homotopic to a closed loop positively transversal to the foliation $\hat{\mathcal{W}}$, and crossing the leafs $\hat{\mathcal{W}}_{T_l(\hat x_i)}$ for $1\le l \le 3$. By the choice of $\varepsilon_i$, for each $l$, $\beta$ has a nonzero winding number with either any point in $\alpha( \hat{\mathcal{W}}_{T_l(\hat x_i)})$ or with any point in $\omega( \hat{\mathcal{W}}_{T_l(\hat x_i)})$. We assume, with no loss in generality, that  for any $\hat z \in \alpha( \hat{\mathcal{W}}_{T_1(\hat x_i)})\cup \alpha( \hat{\mathcal{W}}_{T_2(\hat x_i)}),\, W(\beta, \hat z)\not= 0$. The result follows by taking  $T=T_2\circ T_1^{-1}$ and $\hat p$ either a singularity in $\alpha( \hat{\mathcal{W}}_{T_1(\hat x_i)})$ or, in case $\alpha( \hat{\mathcal{W}}_{T_1(\hat x_i)})$ is a simple closed curve, a singularity in its interior
\endproof

Let now $\hat y$ be a $\hat g$ $k$--periodic point whose orbit has diameter greater than $M$ and consider $\hat p,\,T$ given by the previous lemma. Let $h=g^k,\,\hat h=\hat g^k$ and $I_h$ be the isotopy given by concatenating $I_g$ $k$ times.  Then $h(\hat y)$ is fixed, but $\gamma_{\hat y}^{\hat h}$ is not homotopically trivial in $\hat S\setminus \hat X$, and 
the orbit of $\hat y$ by $\hat h$ has nonzero index with both $\hat p$ and $T(\hat p)$. Let $\til S= \hat S/ \langle T \rangle$ which is homeomorphic to the open annulus $\A$, $\pi_1:\hat S\to\til S$ the projection, and $\til h:\til S\to \til S$ be the induced dynamics satisfying $\til h\pi_1=\pi_1\hat h$.  Let $\til p=\pi_1(\hat p)=\pi_1(T(\hat p)),\, \til y =\pi_1(\hat y), \, \til Y=\{\til p, \til y\}$ and $\hat Y=\pi^{-1}(\til Y)$.


\begin{lemma}\label{lm:pseudoanosov}
 $\til h$ is homotopic in $\til S\setminus \til Y$ to a homeomorphism pseudo-Anosov relative to $\til Y$
\end{lemma}

\proof
Let $\gamma$ be a simple closed curve in $\til M=\til S\setminus\til Y$ which is neither homotopic to a topological end of $\til M$, nor homotopically trivial. Note that if $\gamma'$ is another curve which is also neither homotopic to a topological end of $\til M$, nor homotopically trivial, and such that $[\gamma']$ is disjoint from $[\gamma]$, then $\gamma'$ is homotopic to $\gamma$ in $\til M$. Therefore, if $\til f$ is a reducible homeomorphism and $\gamma$ is one of the reducing curves, $\til f$ must leave $\gamma$ invariant.

Let $\til f$ be homotopic to $\til h$ relative to $\til X$. We show $\til f$ cannot leave $\gamma$ invariant which in turn implies,  by Thurston's  classification, the stated lemma. Let $\til J$ be an isotopy between $\til h$ and $\til f$ leaving $\til Y$ fixed and let $I_{\til f}=I_{\til g}*\til J$ be the isotopy between the identity and $\til f$, and let $I_{\hat f}$ be the lifted isotopy to $\hat S$. Then both $\hat y$ and $\hat p$ are in $\fix(\hat f)$ and furthermore $I_{\hat f}(\hat p, t)=\hat p,\, t\in [0,1].$ Note that, since $\hat S$ is the universal covering of $\til S,\, \hat f$ commutes with $T$, but not necessarily with every transformation in Deck($S$). Finally, note that the closed loops $\gamma_{\hat y}^{\hat f}$ and $\gamma_{\hat y}^{\hat h}$ are homotopic.

\begin{figure}[ht]
\begin{center}\includegraphics[height=6cm]{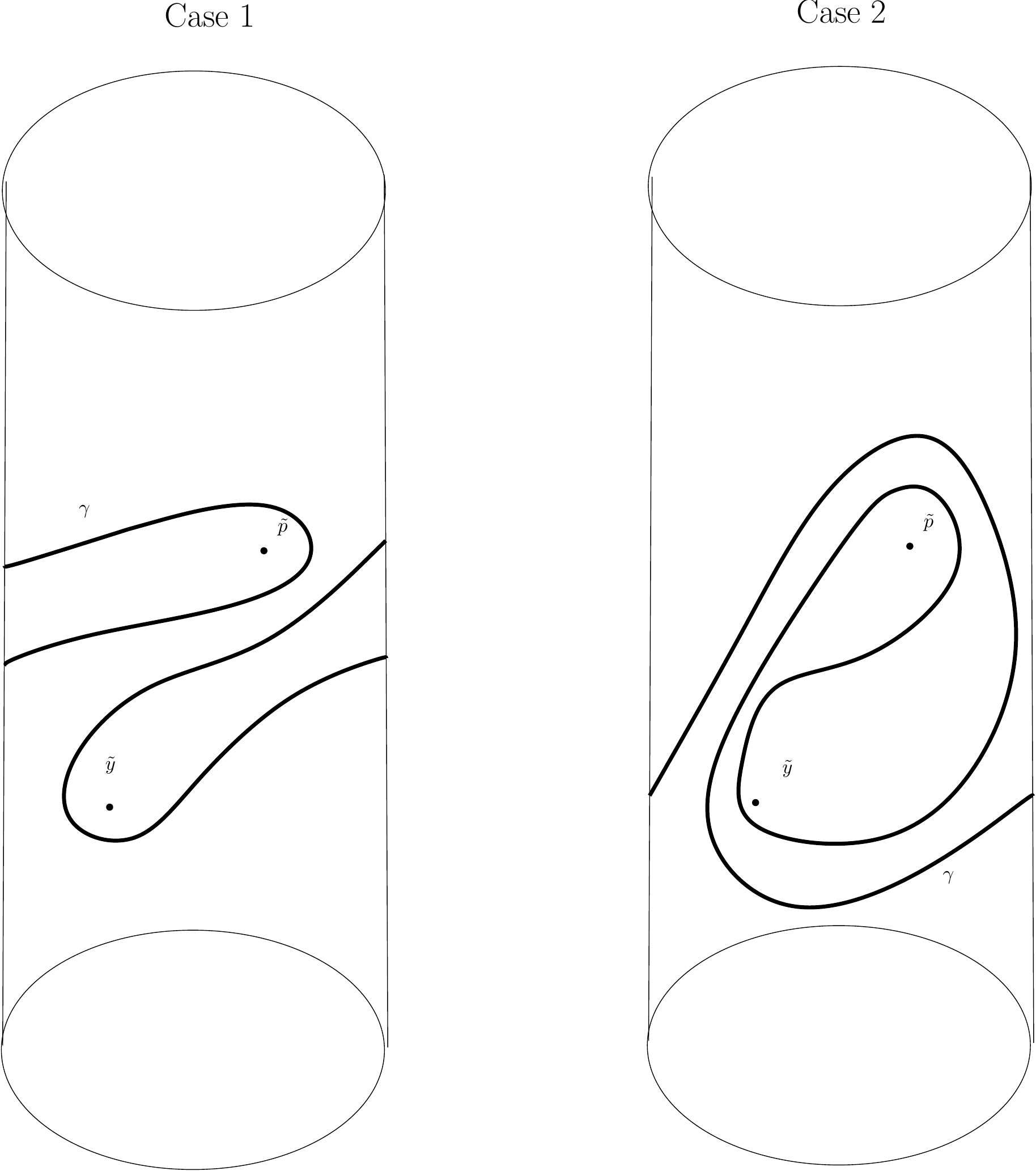}
\end{center}
\caption{Possible reducing curves}
\label{fig:1}
\end{figure}

There are two possibilities to consider
\begin{enumerate}
\item {$\gamma$ is homotopically trivial in $\til S$, in which case $\til S\setminus \gamma$ has two connected components, one which contains $\til p$ and $\til y$ and another which contains the ends of $\til S$,}
\item{$\gamma$ is not homotopically trivial, in which case it separates the ends of $\til S$, leaving $\til p$ in one side and $\til y$ in the other,}
\end{enumerate} 
see figure \ref{fig:1}.

First assume (1). Let $[\hat \gamma]$ be the connected component of  $\pi^{-1}([\gamma])$ such that $\hat y$ is in its interior.
Since $\gamma$ is homotopically trivial in $\til S$, $\hat \gamma$ is also homotopically trivial in $\hat S$, and its interior can contain either $\hat p$ or $T(\hat p)$, but not both. Since 
$\beta(t)=\gamma_{\hat y}^{\hat h}=\gamma_{\hat y}^{\hat f}$ and since we now that $\beta$  has nonzero winding number around  both $\hat p$ and $T(\hat p)$, $\hat f$ cannot keep $\hat \gamma$ invariant.

Now assume (2). In this case, $\hat y$ and $\hat p$ are in different unbounded connected components of $\hat S\setminus \pi_1^{-1}(\gamma)$, and again, since the winding number of $\beta$ around $\hat p$ is not zero, $\hat f$ also cannot leave $\hat \gamma$ invariant.
\endproof

Now, since both $\hat p$ and $\hat y$ are fixed by $\hat h$, and $\hat h$ is the lift of $\til h$ to the universal covering of $\til S$, then both $\rho(\til p, \hat h)$ and $\rho(\til y, \hat h)$ are null. Therefore, Lemmas \ref{lm:toby2} and Lemma \ref{lm:pseudoanosov} imply that, if $\abs{\frac{k p}{q}}$ is sufficiently small, then there exists $\hat z_1 \in \hat S$ such that $\hat g^{k q}(\hat z_1)=\hat h^q(\hat z_1)= T^{k p}(\hat z_1).$ Therefore, as $\hat g$ and $T$ commute, $\hat z_1$ is a $k$-periodic point for $ T^{-p}(\hat g^k)$ and by Brouwer theory there exists $\hat z_2$ fixed by this same map, which ends the proof of Theorem \ref{th:pontos periodicos}. 

\




\subsection{ Proof of Theorem \ref{th:pontos inessenciais}}

Assume that $\pi(\fix(\hat g))$ is inessential. Let $\mathcal{W}$ be a Brouwer - Le Calvez foliation for $\hat g$, and let $ \{x_1\, ...,\,x_n\}, M$ be as in proposition \ref{pr:visita3}.

Assume that there exists $\hat z$ be such that $z=\pi(\hat z)\in\ine(g)\cap\Omega(g)$, and such that the diameter of the $g$ orbit of $\hat z$ is larger than $2M$, otherwise we are done.
Since $z$ is an inessential nonwandering point, there exists a topological disk $U\subset S$ and a integer $k$ such that $g^{k}(U)=U$, and $z\in U$. Since $U$ has a nonwandering point, there must also exists a $k$--periodic point $q$ in $U$.

Let $\hat U$ be a connected component of $\pi^{-1}(U)$. Then there exists a deck transformation $T$ such that $\hat g^k(\hat U)=T(U)$, and in particular, $\hat g^k(\hat q)=T(\hat q)$. If $T$ is not the identity, we are done, so assume that $\hat g^{k}(\hat U)=\hat U$. Since $\hat z\in \hat U$ and $\pi(\hat z)$ is nonwandering, $\hat z$ itself is nonwandering. Let $\delta>0$ be sufficiently small such that, for every $\hat x\in B=B_{\delta}(\hat z)$ the orbit of $\hat x$ has diameter greater than $M$, and that $B\subset \hat U$. We also assume, by proposition \ref{pr:bolaefolhas}, that $\delta$ is sufficiently small such that, if $h$ is a homeomorphism whose restriction to $S\setminus  g^{-1}(B)$ is equal to $g$, then there exists $I_h$ an isotopy between $h$ and the identity, and such that $\mathcal{W}$ is a Le Calvez foliation for $(I_h, X)$.

Let $l=\min\{j>0\mid \hat g^{j}(B)\cap B\not=\emptyset\}$, and let $\hat y \in B$ such that $\hat g^{l}(\hat y)\in B$. Let $h$ be isotopic to $g$, such that $h(x)=g(x)$ for all $x\notin g^{-1}(B)$ and such that $h(g^{l-1}(y))= h^{l}(y)=y$, and let $\hat h$ be the associated lift of $h$. Since $\hat y$ is now a $l$ periodic point for $\hat h$, and since the isotopy path of $\hat g$ between $y$ and $\hat g^{l-1}(y)$ is equal to the isotopy path of $h$ between $\hat y$ and $\hat h^{l-1}(\hat y)$, The exact same argument done in lemma \ref{lm:enlaca_p_e_tp} and the previous theorem yields that there exists $T_1\in Deck(S),\, T_1\not= Id,\, \hat x_1$ and $i$ such that $\hat h^i(\hat x_1)=T(\hat x_1)$. As the orbit of $\hat U$ does not intersect $T_1(\hat U)$, the orbit of $\pi(\hat x_1)$ cannot enter $U$. As $h$ only differs from $g$ on $g^{-1}(U)$, the $h$ orbit and the $g$ orbit of $x_1$ are equal, and so $\hat g^{i}(\hat x_1)=T_1(\hat x_1)$  \endproof

\subsection{Proof of Theorem \ref{th:irrotacional}}

Assume $\fix(g)$ is inessential. In particular, due to Theorem \ref{th:pseudorotacoes}, either $g$  if $g$ is not strictly toral, and thus $g$ is annular, in which case we are done, or every $\hat g$ orbit is bounded, but for any $D>0$ there exists an orbit with diameter larger then $D$. Assume for a contradiction that $g$ is not annular. Since $g$ is area preserving, $\Omega(g)=\T^2$ and so Theorem \ref{th:strictlytoral} applies. 

Also, since $g$ is irrotational, every $g$\--periodic point lifts to a $\hat g$\--periodic point. Theorem \ref{th:pontos inessenciais} implies that, there exists $M>0$ such that, if $\pi(\hat x)\in\ine(g)$, then the orbit of $\hat x$ has diameter less then $M$. Since the orbits by $\hat g$ are not uniformly bounded, this implies that $\cl(\ine(g))$ is not $\T^2$, and thus $\ess(g)$ has nonempty interior. Again by Theorem \ref{th:strictlytoral}, if $U, V$ are open subsets of $\inter(\ess(g))$, then there exists $n>0$ such that $\hat g^n(U)\cap V\not=\emptyset$, so the restriction of $g$ to $\inter(\ess(g))$ is transitive, and there exists $y$ with a dense orbit in $\inter(\ess(g))$. But the orbit of $y$ must have bounded diameter $D_1$, and so the orbit of any point in $\inter(\ess(g))$ also has bounded diameter $D_1$. Since $\T^2=\inter(\ess(g))\cup \cl(\ine(g))$, this implies that the orbit of any point has diameter smaller than $\max\{M, D_1\}$, a contradiction.

\section{Aknowledgements}
 We would like to thank S. Addas--Zanata, T. Hall, and A. Koropecki for the insightful discussions that helped simplify this work.


\bibliographystyle{amsalpha}

\end{document}